\documentstyle[12pt]{amsart}
\newcommand{\g}{\goth}
\newcommand{\gtg}{\mbox{\g g}}

\newcommand{\hgtg}{\mbox{$\hat{\gtg}$}}

\newcommand{\gtsl}{\mbox{\g sl}}

\newcommand{\hgtsl}{\mbox{$\hat{\gtsl}$}}

\newcommand{\nc}{\mbox{${\Bbb C}$}}
\newcommand{\nz}{\mbox{${\Bbb Z}$}}

\newcommand{\cF}{\mbox{${\cal F}$}}

\newcommand{\cP}{\mbox{${\cal P}$}}

\newcommand{\vep}{\varepsilon}
\newcommand{\vphi}{\varphi}
\theoremstyle{plain}
 \newtheorem{thm}{Theorem}[section]
 \newtheorem{prop}[thm]{Proposition}
 \newtheorem{lemma}[thm]{Lemma}
 
\theoremstyle{definition}
 \newtheorem{defn}{Definition}[section]
 
 \newtheorem{ex}{Example}[section]
\theoremstyle{remark}

 \newtheorem{ack}{Acknowledgment}
      
\begin{document}
\title{Generalization and Deformation of Drinfeld quantum affine algebras}
\author{Jintai Ding}
\author{kenji Iohara}
 
\address{Jintai Ding, RIMS, Kyoto University}
\address{Kenji Iohara, Dept. of Math., Kyoto University}

\maketitle
\begin{abstract}
Drinfeld gave a current realization of the quantum affine algebras as a 
Hopf algebra   with a simple
 comultiplication for the quantum current operators. 
In this paper, we will present a generalization of such a realization 
of quantum Hopf algebras. 
As a special case, we will  choose the structure functions for this 
algebra to be elliptic functions to derive certain  elliptic quantum groups
as a Hopf algebra, which degenerates into quantum affine algebras if we take 
certain degeneration of the structure functions. 
\end{abstract}

\pagestyle{plain}
\section{Introduction.} 

Quantum groups  as the very first noncommutative 
and noncocommutative  Hopf algebras were   discovered by Drinfeld\cite{Dr1} 
and Jimbo\cite{J1}. 
The standard  definition of a  quantum group is given  as a deformation 
of a simple Lie algebra
 by the basic generators and the relations based on the data coming from the  
corresponding Cartan matrix. However, for the case of affine quantum 
groups, there is a different aspect of the theory. 

For the  undeformed  affine Lie algebras,
Garland showed  \cite{G} that the affine
Kac-Moody algebra $\hat {\frak{g}}$ associated to a simple Lie algebra
$\frak{g}$ admits a natural realization as a central extension of the
corresponding loop algebra ${\frak{ g}} \otimes \Bbb C [t, t^{-1}] $.
It is natural to expect similar realizations  for the quantum affine algebras. 
The first approach was given by Faddeev, Reshetikhin and
Takhtajan \cite{FRT}, who obtained a realization  of
the quantum loop algebra $U_q(
{\frak{ g}} \otimes \Bbb [t, t^{-1}])$ via a canonical 
 solution
 of the Yang-Baxter equation depending on a parameter $z\in\Bbb C$. 
This approach was completed by Reshetikhin and Semenov-Tian-Shansky
\cite{RS} by  incorporating  the central extension in the
previous realization. However, in this approach, we can only have 
either the positive current operators or the negative current operators, 
which have complicated commutation relations. 

 On the other hand, Drinfeld \cite{Dr2} gave another realization  
of the quantum affine algebra $U_q(\hgtg)$ and its special degeneration 
called the Yangian, 
which immediately is used in constructions of special representations of
affine quantum algebras \cite{FJ}.  
In \cite{Dr2}, Drinfeld only gave the realization of 
the quantum affine algebras as an algebra and as an  algebra this 
realization is equivalent to  the approach above \cite{DF} through
certain Gauss decomposition for the case of $U_q(\hat {\frak gl}(n))$. 
Certainly, the most important aspect of the structures 
 of the quantum groups is its Hopf algebra
structure, especially its comultiplication. 
However, if we  extend the conventional comultiplication to 
the current operators,
we acquire a very complicated formula which  can not be written 
in a closed form with only the current operators. 
Drinfeld also gave the Hopf algebra structure 
for such a formulation in an unpublished note \cite{DF}.
The new  comultiplication in this formulation, which we call  the  Drinfeld 
comultiplication, is simple. 
Recently in \cite{DM} \cite{DI}, we used this new Hopf algebra structure
to study vertex operators and zeros and poles of the quantum current 
operators. 

In the Drinfeld realization  of 
quantum affine algebras, the structure constants are certain 
rational functions $g_{ij}(z)$. 
In  checking the Hopf algebra structure  in this formulation, 
we notice that   centain functional property of $g_{ij}(z)$ decides that 
this  formulation gives 
us a Hopf algebra. This leads  us to 
generalize this type of Hopf algebras. Namely, we can substitute  
$g_{ij}(z)$ by other functions that satisfy the functional property 
of $g_{ij}(z)$, to derive some  new Hopf algebras. 

In this paper, we will present our generalization of
the Drinfeld realization of  $U_q(\hgtsl_n)$. 
We will show that   $U_q(\hgtsl_n)$
are simple examples of such a formulation. 
We will then choose the structure functions for this 
algebra to be certain  elliptic functions.
We  derive an elliptic quantum groups
as a Hopf algebra, which degenerates into quantum affine algebras if we choose
special degeneration of the structure functions.

\section{}
The Drinfeld realization for the case of 
$U_q(\hgtsl_n)$ \cite{Dr2} as a Hopf algebra is presented in \cite{DF}
 \cite{DI}, 
which leads  us to define the following algebra. 

Let $A=(a_{ij})$ be the Cartan matrix of type $A_{n-1}$.
Let $g$ be the set of 
$\{g_{ij}(z),0<i,j<n \}$
which are  analytic functions  satisfying the 
following property that $g_{ij}(z)=g_{ji}(z^{-1})^{-1}$ and
$\delta(z)$ be  the distribution with support at $1$.

\begin{defn}
$U_q(g,\gtsl_n)$ is an
 associative algebra with unit 
1 and the generators:  $x_i^{\pm}(z)$, 
$\vphi_i(z)$,
$\psi_i(z)$, a central element $c$
 and a nonzero  complex parameter $q$. 
$\vphi_i(z)$ and $\psi_i(z)$ are invertible. 
In terms of the generating functions:
the defining relations are 
\begin{align*}
& \vphi_i(z)\vphi_j(w)=\vphi_j(w)\vphi_i(z), \\
& \psi_i(z)\psi_j(w)=\psi_j(w)\psi_i(z), \\
& \vphi_i(z)\psi_j(w)\vphi_i(z)^{-1}\psi_j(w)^{-1}=
  \frac{g_{ij}(z/wq^{-c})}{g_{ij}(z/wq^{c})}, \\
& \vphi_i(z)x_j^{\pm}(w)\vphi_i(z)^{-1}=
  g_{ij}(z/wq^{\mp \frac{1}{2}c})^{\pm1}x_j^{\pm}(w), \\
& \psi_i(z)x_j^{\pm}(w)\psi_i(z)^{-1}=
  g_{ij}(w/zq^{\mp \frac{1}{2}c})^{\mp1}x_j^{\pm}(w), \\
& [x_i^+(z),x_j^-(w)]=\frac{\delta_{i,j}}{q-q^{-1}}
  \left\{ \delta(z/wq^{-c})\psi_i(wq^{\frac{1}{2}c})-
          \delta(z/wq^{c})\vphi_i(zq^{\frac{1}{2}c}) \right\}, \\
&  x_i^{\pm}(z)x_j^{\pm}(w)=
  g_{ij}(z/w)^{\pm 1}x_j^{\pm}(w)x_i^{\pm}(z), \\
& [x_i^{\pm}(z),x_j^{\pm}(w)]=0 \quad \text{ for $a_{ij}=0$}, \\
& x_i^{\pm}(z_1)x_i^{\pm}(z_2)x_j^{\pm}(w)-(h_{ij}(z_1/w,z_2/w))x_i^{\pm}(z_1)
  x_j^{\pm}(w)x_i^{\pm}(z_2)+x_j^{\pm}(w)x_i^{\pm}(z_1)x_i^{\pm}(z_2) \\
& +\{ z_1\leftrightarrow z_2\}=0, \quad \text{for $a_{ij}=-1$}
\end{align*}
where
$$  h_{ij}(z_1/w,z_2/w)=\frac {(g_{ii} (z_1/z_2)+1)(g_{ij}(z_1/w)g_{ij}(z_2/w)
+1)}
{g_{ij}(z_2/w)+g_{ii}(z_1/z_2)g_{ij}(z_1/w)} 
 \]
\end{defn}

\begin{thm}
The algebra $U_q(g,\gtsl_n)$  
has a Hopf algebra structure, which are given 
by the following formulae. 

\noindent{\bf Coproduct $\Delta$}
\begin{align*}
\text{(0)}& \quad \Delta(q^c)=q^c\otimes q^c, \\
\text{(1)}& \quad \Delta(x_i^+(z))=x_i^+(z)\otimes 1+
            \vphi_i(zq^{\frac{c_1}{2}})\otimes x_i^+(zq^{c_1}), \\
\text{(2)}& \quad \Delta(x_i^-(z))=1\otimes x_i^-(z)+
            x_i^-(zq^{c_2})\otimes \psi_i(zq^{\frac{c_2}{2}}), \\
\text{(3)}& \quad \Delta(\vphi_i(z))=
            \vphi_i(zq^{-\frac{c_2}{2}})\otimes\vphi_i(zq^{\frac{c_1}{2}}), \\
\text{(4)}& \quad \Delta(\psi_i(z))=
            \psi_i(zq^{\frac{c_2}{2}})\otimes\psi_i(zq^{-\frac{c_1}{2}}),
\end{align*}
where $c_1=c\otimes 1$ and $c_2=1\otimes c$.

\noindent{\bf Counit $\vep$}
\begin{align*}
\vep(q^c)=1 & \quad \vep(\vphi_i(z))=\vep(\psi_i(z))=1, \\
            & \quad \vep(x_i^{\pm}(z))=0.
\end{align*}
\noindent{\bf Antipode $\quad a$}
\begin{align*}
\text{(0)}& \quad a(q^c)=q^{-c}, \\
\text{(1)}& \quad a(x_i^+(z))=-\vphi_i(zq^{-\frac{c}{2}})^{-1}
                               x_i^+(zq^{-c}), \\
\text{(2)}& \quad a(x_i^-(z))=-x_i^-(zq^{-c})
                               \psi_i(zq^{-\frac{c}{2}})^{-1}, \\
\text{(3)}& \quad a(\vphi_i(z))=\vphi_i(z)^{-1}, \\
\text{(4)}& \quad a(\psi_i(z))=\psi_i(z)^{-1}.
\end{align*}

\end{thm}

 \noindent{Proof.}
For the comultiplication above we have that 
\begin{align*}
\text{(1)}& \quad \Delta(x_i^+(z))=x_i^+(z)\otimes 1+
            \vphi_i(zq^{\frac{c_1}{2}})\otimes x_i^+(zq^{c_1}), \\
\text{(2)}& \quad \Delta(x_i^-(z))=1\otimes x_i^-(z)+
            x_i^-(zq^{c_2})\otimes \psi_i(zq^{\frac{c_2}{2}}), \\
\text{(3)}& \quad \Delta(\vphi_i(z))=
            \vphi_i(zq^{-\frac{c_2}{2}})\otimes\vphi_i(zq^{\frac{c_1}{2}}), \\
\text{(4)}& \quad \Delta(\psi_i(z))=
            \psi_i(zq^{\frac{c_2}{2}})\otimes\psi_i(zq^{-\frac{c_1}{2}}).
\end{align*}
 $$\Delta\vphi_i(z)\Delta\vphi_j(w)=\Delta\vphi_j(w)\Delta\vphi_i(z),$$
$$\Delta\psi_i(z)\Delta\psi_j(w)=\Delta\psi_j(w)\Delta\psi_i(z).$$
Then, 
$$ \Delta\vphi_i(z)\Delta\psi_j(w)\Delta\vphi_i(z)^{-1}\Delta\psi_j(w)^{-1}=$$
  $$\vphi_i(zq^{-\frac{c_2}{2}})\psi_j(wq^{\frac{c_2}{2}})
\vphi_i(zq^{-\frac{c_2}{2}})^{-1}\psi_j(wq^{\frac{c_2}{2}})^{-1}
\otimes  \vphi_i(zq^{\frac{c_1}{2}})\psi_j(wq^{-\frac{c_1}{2}})
\vphi_i(zq^{\frac{c_1}{2}})^{-1}\psi_j(wq^{-\frac{c_1}{2}})^{-1}=$$
$$   \frac{g_{11}(z/wq^{-c_1}q^{-{c_2}} )}
{g_{{ij}}(z/wq^{c_1}q^{-{c_2}})}
\frac{g_{{ij}}(z/wq^{-c_2}q^{{c_1}})}{g_{{ij}}(z/wq^{c_2}q^{{c_1}})}=$$
$$\frac{g_{{ij}}(z/wq^{-c_1-c_2} )}
{g_{{ij}}(z/wq^{c_1+c_2})}.$$

$$  \Delta\vphi_i(z)\Delta x_j^{+}(w)\Delta\vphi_i(z)^{-1}=$$
$$\vphi_i(zq^{-\frac{c_2}{2}})\otimes\vphi_i(zq^{\frac{c_1}{2}})
(x_j^+(w)\otimes 1+
            \vphi_j(wq^{\frac{c_1}{2}})\otimes x_j^+(wq^{c_1}))
 (\vphi_i(zq^{-\frac{c_2}{2}})\otimes\vphi_i(zq^{\frac{c_1}{2}}))^{-1}=$$
$$ g_{ij}(z/wq^{ -\frac{1}{2}(c_1+c_2)})(x_j^+(w)\otimes 1+
            \vphi_j(wq^{\frac{c_1}{2}})\otimes x_j^+(wq^{c_1})). $$

$$  \Delta\vphi_i(z)\Delta x_j^{-}(w)\Delta\vphi_i(z)^{-1}=$$
$$\vphi_i(zq^{-\frac{c_2}{2}})\otimes\vphi_i(zq^{\frac{c_1}{2}})
(1\otimes x_j^-(w)+
            x_j^-(wq^{c_2})\otimes \psi_j(wq^{\frac{c_2}{2}}))
(\vphi_i(zq^{-\frac{c_2}{2}})\otimes\vphi_i(zq^{\frac{c_1}{2}}))^{-1}=$$
$$1\otimes x_j^-(w)g_{ij}(z/wq^{ \frac{1}{2}(c_1+c_2)})^{-1}
+x_j^-(wq^{c_2})\otimes \psi_j(wq^{\frac{c_2}{2}})
g_{ij}(z/wq^{ \frac{1}{2}(c_1-3c_2)})^{-1}
\frac {g_{ij}(z/wq^{ \frac{1}{2}(c_1-3c_2)})}
      {g_{ij}(z/wq^{ \frac{1}{2}(c_1+c_2)})}=$$
$$\Delta x_i^{-}(w)g_{ij}(z/wq^{ \frac{1}{2}(c_1+c_2)})^{-1}.$$

The relation between $\Delta\psi_i(z)$ and $\Delta x_i^{\pm}(w)$
 can be proved with the same method demonstrted above. 

$$ [\Delta x_i^+(z),\Delta x_i^-(w)]=$$
$$[x_i^+(z)\otimes 1+
            \vphi_i(zq^{\frac{c_1}{2}})\otimes x_i^+(zq^{c_1}),$$
$$ 
1\otimes x_i^-(w)+
            x_i^-(wq^{c_2})\otimes \psi_i(wq^{\frac{c_2}{2}})]=$$
$$[x_i^+(z)\otimes 1, x_i^-(wq^{c_2})\otimes \psi_i(wq^{\frac{c_2}{2}})]+
[\vphi_i(zq^{\frac{c_1}{2}})\otimes x_i^+(zq^{c_1}), 
1\otimes x_i^-(w)]+$$
$$ [\vphi_i(zq^{\frac{c_1}{2}})\otimes x_i^+(zq^{c_1}), 
            x_i^-(wq^{c_2})\otimes \psi_i(wq^{\frac{c_2}{2}})].$$
It is easy to show that the last term above is 0, therefore we have that
$$ [\Delta x_i^+(z),\Delta x_i^-(w)]=$$
$$\frac{(\delta(z/wq^{-c_1-c_2})\psi_i(wq^{\frac{c_1}{2}+c_2})-
\delta(z/wq^{c_1-c_2})\vphi_i(zq^{\frac{c_1}{2}}))\otimes 
\psi_i(wq^{\frac{c_2}{2}})}{q-q^{-1}}+$$
$$\frac{\vphi_i(zq^{\frac{c_1}{2}})\otimes
 ( \delta(z/wq^{c_1-c_2})\psi_i(wq^{\frac{1}{2}c_2})-
          \delta(z/wq^{c_1+c_2})\vphi_i(zq^{c_1+\frac{1}{2}c_2}))}{q-q^{-1}}=$$
$$ \frac{ \delta(z/wq^{-c_1-c_2})\psi_i(zq^{\frac{c_2}{2}+\frac{c_1+c_2}{2}})
\otimes\psi_i(zq^{-\frac{c_1}{2}+\frac{c_1+c_2}{2}})-
\delta(z/wq^{c_1+c_2})\vphi_i(zq^{-\frac{c_2}{2}+\frac{c_1+c_2}{2}})
\otimes\vphi_i(zq^{\frac{c_1}{2}+\frac{c_1+c_2}{2}})}{q-q^{-1}}.$$

$$\Delta x_i^{+}(z)\Delta x_j^{+}(w)=$$
$$(x_i^+(z)\otimes 1+
            \vphi_i(zq^{\frac{c_1}{2}})\otimes x_i^+(zq^{c_1}))
(x_j^+(w)\otimes 1+
            \vphi_j(wq^{\frac{c_1}{2}})\otimes x_j^+(wq^{c_1}))=$$ 
$$ g_{ij}(z/w)x_j ^{+}(w)x_i^{+}(z)\otimes 1+
\vphi_j(wq^{\frac{c_1}{2}}x_i^+(z)\otimes x_j^+(wq^{c_1})
g_{ji}(w/z)^{-1})+$$
$$
x_j^+(w)\vphi_i(zq^{\frac{c_1}{2}})\otimes  
x_i^+(zq^{c_1})g_{ij}(z/w)+$$
$$ g_{ij}(z/w)\vphi_j(wq^{\frac{c_1}{2}})\vphi_i(zq^{\frac{c_1}{2}})\otimes 
   x_j^+(wq^{c_1}) x_i^+(zq^{c_1})=$$
$$g_{ij}(z/w)\Delta x_i^{+}(w)\Delta x_j^{+}(z).$$

The last is to check the cubic relations between 
$\Delta x^{\pm}_i(z)$ and $\Delta x^{\pm}_j(w)$ for $a_{ij}=-1$. 
Therefore  we need to show that 
$$ \Delta x_i^{+}(z_1)\Delta x_i^{+}(z_2)\Delta x_j^{+}(w)
-(h_{ij}(z_1/w,z_2/w))\Delta x_i^{+}(z_1)\Delta x_j^{+}(w)
\Delta x_i^{+}(z_2)+$$ 
$$\Delta x_j^{+}(w)\Delta x_i^{+}(z_1)\Delta x_i^{+}(z_2)
+\{ z_1\leftrightarrow z_2\}=0 $$

$$(x_i^+(z_1)\otimes 1+
            \vphi_i(z_1q^{\frac{c_1}{2}})\otimes x_i^+(z_1q^{c_1}))
(x_i^+(z_2)\otimes 1+
            \vphi_i(z_2q^{\frac{c_1}{2}})\otimes x_i^+(z_2q^{c_1}))\times $$
$$
(x_j^+(w)\otimes 1+
            \vphi_j(wq^{\frac{c_1}{2}})\otimes x_j^+(wq^{c_1}))=$$
$$ x_i^+(z_1)x_i^+(z_2)x_j^+(w)\otimes 1+ x_i^+(z_1)x_i^+(z_2)
\vphi_j(wq^{\frac{c_1}{2}})\otimes x_j^+(wq^{c_1})+$$
$$x_i^+(z_1) \vphi_i(z_2q^{\frac{c_1}{2}})x_j^+(w)\otimes 
 x_i^+(z_2q^{c_1})+x_i^+(z_1) \vphi_i(z_2q^{\frac{c_1}{2}}) 
\vphi_j(wq^{\frac{c_1}{2}})\otimes x_i^+(z_2q^{c_1})x_j^+(wq^{c_1})+$$
$$ \vphi_i(z_1q^{\frac{c_1}{2}})x_i^+(z_2)x_j^+(w)\otimes x_i^+(z_1q^{c_1})+
 \vphi_i(z_1q^{\frac{c_1}{2}})x_i^+(z_2)
\vphi_j(wq^{\frac{c_1}{2}})\otimes x_i^+(z_1q^{c_1}) x_j^+(wq^{c_1})+$$
$$\vphi_i(z_1q^{\frac{c_1}{2}}) \vphi_i(z_2q^{\frac{c_1}{2}})x_j^+(w)\otimes
 x_i^+(z_1q^{c_1}) 
 x_i^+(z_2q^{c_1})+$$
$$\vphi_i(z_1q^{\frac{c_1}{2}}) \vphi_i(z_2q^{\frac{c_1}{2}}) 
\vphi_j(wq^{\frac{c_1}{2}})\otimes x_i^+
(z_1q^{c_1}) x_i^+(z_2q^{c_1})x_j^+(wq^{c_1}).$$

Note that $$h_{ij}(z_1/w,z_2/w)=h_{ij}(z_2/w,z_1/w).$$

Thus we have that
\begin{align*}
&  \Delta x_i^{+}(z_1)\Delta x_i^{+}(z_2)\Delta x_j^{+}(w)
-(h_{ij}(z_1/w,z_2/w))\Delta x_i^{+}(z_1)\Delta x_j^{+}(w)
\Delta x_i^{+}(z_2)+ \\
&  \Delta x_j^{+}(w)\Delta x_i^{+}(z_1)\Delta x_i^{+}(z_2) 
+\{ z_1\leftrightarrow z_2\}= \\
& \left\{ 
\vphi_j(wq^{\frac{c_1}{2}})x_i^+(z_2)x_i^+(z_1)\otimes x_j^+(wq^{c_1})+
x_j^+(w)x_i^+(z_2)\vphi_i(z_1q^{\frac{c_1}{2}})\otimes x_i^+(z_1q^{c_1})+
  \right. \\
& x_j^+(w)\vphi_i(z_2q^{\frac{c_1}{2}})x_i^{+}(z_1)\otimes x_i^+(z_2q^{c_1})+
x_j^+(w)\vphi_i(z_2q^{\frac{c_1}{2}})\vphi_i(z_1q^{\frac{c_1}{2}})\otimes
  x_i^+(z_2q^{c_1})x_i^+(z_1q^{c_1}) \\
&  
\vphi_j(wq^{\frac{c_1}{2}})x_i^+(z_2)\vphi_i(z_1q^{\frac{c_1}{2}})\otimes
x_j^+(wq^{c_1})x_i^+(z_1q^{c_1})+ \\
& \left.
\vphi_j(wq^{\frac{c_1}{2}})\vphi_i(z_2q^{\frac{c_1}{2}})x_i^{+}(z_1)\otimes
x_j^+(wq^{c_1})\otimes x_i^+(z_2q^{c_1}) \right\}\times \\
& \left\{
g_{ii}(z_1/z_2)  g_{ij}(z_1/w)g_{ij}(z_2/w)-
h_{ij}(z_1/w,z_2/w)g_{ii}(z_1/z_2)g_{ij}(z_1/w)+g_{ii}(z_1/z_2)+ \right. \\
& \left.
g_{ij}(z_1/w)g_{ij}(z_2/w)-h_{ij}(z_2/w, z_1/w)g_{ij}(z_2/w)+1 \right\} 
\end{align*}

We therefore
 derive the proof that the comultiplication is an algebra homomorphism. 

We can prove in the same manner   the relation between 
$\Delta x_i^-(z)$ and $\Delta x_j^-(w)$.

Let $M$ be the operator from $U_q(g,\gtsl_n)\otimes U_q(g,\gtsl_n)$ to 
$U_q(g,\gtsl_n)$ defined by the algebra multiplication. We can check that
$M(1\otimes \vep)\Delta=\text{ id}$. 

$$ M(1\otimes a) \Delta(x_i^+(z))=$$
$$M(1\otimes a)( x_i^+(z)\otimes 1+
            \vphi_i(zq^{\frac{c_1}{2}})\otimes x_i^+(zq^{c_1}))=$$
$$ x_i^+(z)-\vphi_i(zq^{\frac{c}{2}})(\vphi_i(zq^{\frac{c}{2}}))^{-1}
 x_i^+(z)=$$
$$0=\vep(x_i^+(z))$$
Similarly we can check all the other relations to show that 
$$M(a\otimes 1)=\vep.$$ 
Thus, we prove that the comultiplication, the counit and the 
antipode give a Hopf algebra
structure. 
 
Strictly speaking, $U_q(g,\gtsl_n)$ is not an algebra. 
This concept, which  we call a functional algebra, has  already 
been used beofre \cite{S}, etc. 
For our algebra $U_q(g,\gtsl_n)$,
each generating functions 
gives an operator at any point $\Bbb C$, or if we 
restrict our generating functions at any  point, we get an associative
algebra. Certainly, this is not the case at the poles of 
the structure functions $g_{ij}(z)$. A functional algebra as
a generalization of algebra is similar to the generalization from the concept 
of  a number to a function. Furthermore,  we can
 extend the definition of 
this functional algebra to any abelian 
Lie group. We   simply  substitute 
$\Bbb C$ by an  abelian Lie group, for example:  an elliptic curve, set 
$z$ and  $q^{c/2}$ by two points on the abelian Lie group,
take $g_{ij}(z)$  as functions on the abelian Lie group and 
$\delta(z)$ the distribution with support at point $1$ and 
identify $z^{-1}$ as the inverse of $z$ as a group element. 
Hoever, this should not  be called 
a Hopf algebra, rather a functional Hopf algebra. 

The Drinfeld realization for the case of 
$U_q(\hgtsl_n)$ \cite{Dr2} as a Hopf algebra is different from 
the definition above. It is an algebra and Hopf algebra. 
We will present similar generalization in the same way.

Let $A=(a_{ij})$ be the Cartan matrix of type $A_{n-1}$.
Let $g$ be the set of 
$\{g_{ij}(z),0<i,j<n \}$
be analytic functions that satisfying the 
following property that $g_{ij}(z)=g_{ji}(z^{-1})^{-1}=
G^+_{ij}(z)/G^-_{ij}(z)$, where $G^\pm_{ij}(z)$ is an analytic function
 without poles except at $0$ or $\infty$ and $G^\pm_{ij}(z)$ have no common
zero point.
Let $ \delta(z)=\sum_{n\in \Bbb Z}z^n$.

\begin{defn}
The algebra $U_q(g,\gtsl_n)$ is an associative algebra with unit 
1 and the generators: $\bar a_i(l)$,$\bar b_i(l)$, $x^{\pm}_i(l)$, for 
$i=i,...,n-1$, $l\in \nz $ and a central element $c$. 
Let $z$ be a formal variable and 
$$x_i^{\pm}(z)=\sum_{l\in \nz}x_i^{\pm}(l)z^{-l}$$, 
$$\vphi_i(z)=\sum_{m\in \nz}\vphi_i(m)z^{-m}=\exp [
\sum_{m\in \nz_{\leq 0}}\bar a_i(m)z^{-m}]\exp [
\sum_{m\in \nz_{< 0}}\bar a_i(m)z^{-m}]$$
  and 
$$\psi_i(z)=\sum_{m\in \nz}\psi_i(m)z^{-m}=\exp [
\sum_{m\in \nz_{\leq 0}}\bar b_i(m)z^{-m}]\exp [
\sum_{m\in \nz_{< 0}}\bar b_i(m)z^{-m}] $$.
 In terms of the 
formal variables $z$, $w$, $z_1$ and $z_2$, 
the defining relations are 
\begin{align*}
& a_i(l)a_j(m)=a_j(m)a_i(l), \\
& b_i(l)b_j(m)=b_j(m)b_i(l), \\
& \vphi_i(z)\psi_j(w)\vphi_i(z)^{-1}\psi_j(w)^{-1}=
  \frac{g_{ij}(z/wq^{-c})}{g_{ij}(z/wq^{c})}, \\
& \vphi_i(z)x_j^{\pm}(w)\vphi_i(z)^{-1}=
  g_{ij}(z/wq^{\mp \frac{1}{2}c})^{\pm1}x_j^{\pm}(w), \\
& \psi_i(z)x_j^{\pm}(w)\psi_i(z)^{-1}=
  g_{ij}(w/zq^{\mp \frac{1}{2}c})^{\mp1}x_j^{\pm}(w), \\
& [x_i^+(z),x_j^-(w)]=\frac{\delta_{i,j}}{q-q^{-1}}
  \left\{ \delta(z/wq^{-c})\psi_i(wq^{\frac{1}{2}c})-
          \delta(z/wq^{c})\vphi_i(zq^{\frac{1}{2}c}) \right\}, \\
& G_{ij}^{\mp}(z/w) x_i^{\pm}(z)x_j^{\pm}(w)=
  G_{ij}^{\pm }(z/w)x_j^{\pm}(w)x_i^{\pm}(z), \\
& [x_i^{\pm}(z),x_j^{\pm}(w)]=0 \quad \text{ for $a_{ij}=0$}, \\
& x_i^{\pm}(z_1)x_i^{\pm}(z_2)x_j^{\pm}(w)-(h_{ij}(z_1/w,z_2/w))x_i^{\pm}(z_1)
  x_j^{\pm}(w)x_i^{\pm}(z_2)\\
& +x_j^{\pm}(w)x_i^{\pm}(z_1)x_i^{\pm}(z_2)
  +\{ z_1\leftrightarrow z_2\}=0, \quad \text{for $a_{ij}=-1$}
\end{align*}
where
$$  h_{ij}(z_1/w,z_2/w)=\frac {(g_{ii} (z_1/z_2)+1)(g_{ij}(z_1/w)g_{ij}(z_2/w)
+1)}
{g_{ij}(z_2/w)+g_{ii}(z_1/z_2)g_{ij}(z_1/w)} 
 \]
and  by  $g_{ij}(z)$ we   mean the Laurent expansion of $g_{ij}(z)$ 
in a region $r_{ij,1}>|z|> r_{ij,2}$. 
\end{defn}

\begin{thm}
The algebra $U_q(g,\gtsl_n)$  
has a Hopf algebra structure.  The formulas for the coproduct $\Delta$, 
the counit $\vep$ and the antipode $a$ are the same as given in 
Theorem 2.1. 
\end{thm}

The proof is almost the same as the proof for Theorem 2.1,
except that one has to be careful with the 
expansion direction of the structure functions $g_{ij}(z)$ and 
$\delta(z)$, for the  reason that the  relations between $x^{\pm}_i(z)$
and $x^{\pm}_j(z)$ are different from the case of the  functional algebra
above. 
If we modify the relation between  $x^{\pm}_i(z)$
and $x^{\pm}_j(z)$ in Definiton 2.1 to be the same as in Definition 2.2., 
this will make a subtle difference on the algebra. 
One  problem that appears  in such a formulation is that the  commutation 
relations above can involve infinite expressions, which 
 reqires certain topology to define our algebras. 
 However, it is
 well-defined in the case of the highest weight representations. 

\begin{ex}

Let $g(z)$ be a an analytic function such that $g(z^{-1})=-z^{-1}g(z)$. Let 
$g_{ij}(z)=q^{-a_{ij}}\frac{g(q^{a_{ij}}z)}{g(q^{-a_{ij}}z)}$. Then 
$g_{ij}(z)=g_{ij}(z^{-1})^{-1}$. 
With these $g_{ij}(z)$, we can define
an algebra $U_q(g,\gtsl_n)$ by  Definition 2.2.

For this algebra, we obtain its vector representation. 

Let $V=\oplus_{i=0}^{n-1} \nc|i\rangle$ be a n-dimensional space,
where $\{ |i\rangle \}$ is  its standard basis. Let
$V_z=V\otimes \nc[z,z^{-1}]$, where $z$ is a formal variable.

\begin{lemma}There exists an  n-dimensional representation of 
$U_q(g,\gtsl_n)$
on $V$. The action of the current operators is given by the following:
\begin{align*}
\circ& \quad x_i^+(w).|j\rangle= c_i^+\delta_{ij}
                                 \delta(\frac{w}{q^iz})|i-1\rangle,\
 \\
\circ& \quad x_i^-(w).|j\rangle= c_i^-\delta_{i-1j}
                                 \delta(\frac{w}{q^iz})|i\rangle,\
 \\
\circ& \quad \begin{cases}
             \vphi_i(w).|i-1\rangle
       & =\displaystyle{q^{-1}g(q^2\frac{w}{q^iz})
          \sum_{k\geq 0}(\frac {w}{q^iz})^k}
         |i-1\rangle \\
             \vphi_i(w).|i\rangle
       & =\displaystyle{qg(q^{-2}\frac{w}{q^iz})
          \sum_{k\geq 0}(\frac {w}{q^iz})^k}
         |i\rangle \\
\vphi_i(w).|j\rangle
        & =\displaystyle{g(\frac{w}{q^iz})
          \sum_{k\geq 0}(\frac {w}{q^iz})^k}
         |j\rangle \quad (j\neq i,i-1)
 \end{cases}, \\
\circ& \quad \begin{cases}
             \psi_i(w).|i-1\rangle
       & =\displaystyle{qg(q^{-2}\frac{q^iz}{w})
         \sum_{k\geq 0}(\frac {q^iz}{w})^k}
         |i-1\rangle, \\
             \psi_i(w).|i\rangle
       & =\displaystyle{q^{-1}g(q^2\frac{q^iz}{w})
         \sum_{k\geq 0}(\frac {q^iz}{w})^k}
         |i\rangle \\
 \psi_i(w).|j\rangle
       & =\displaystyle{g(\frac{q^iz}{w})
         \sum_{k\geq 0}(\frac {q^iz}{w})^k}
         |j\rangle \quad (j\neq i,i-1) 
\end{cases}.
\end{align*}
where we set $|-1\rangle =|n\rangle=0$ and $c_i^{\pm}$ are constants satisfying
\[ c_i^+c_i^-=\frac{qg(q^{-2})}{q-q^{-1}}.\]
\end{lemma}

{\bf Case I.} 
Let $g(z)$=$1-z$. 
$\vphi_i(m)=\psi_i(-m)=0$ for $m\in \nz_{>0}$ and $\vphi_i(0)\psi_i(0)=1$. 
Then this 
algebra is  $U_q(\hgtsl_n)$.

In this case, the comultiplication structure requires certain completion
on the tenor space. Nevertheless, if we consider 
any two highest weight representations, which 
can be defined as for the case of $U_q(\hgtsl_n)$, 
this comultiplication
is well-defined. For the case 
of $U_q(\hgtsl_n)$, this 
 is used \cite{DM} to study the poles and zeros
of those current operators for integrable representations. 
In \cite{DI}, this is also used to study the corresponding vertex operators, 
where we derive simple bosonization formulas for vertex operators.

{\bf Case II}
Let $\theta_p(z)=\prod_{j>0}(1-p^j)(1-p^{j-1}z)(1-p^jz^{-1})$ be the Jacobi's
theta function.
$\theta_p(z^{-1})=-z^{-1}\theta_p(z)=\theta_p(pz)$. 
 We will call this 
algebra $U_q(\theta, \gtsl_n)$. We take the expansion of $g_{ii}(z)$
in the region $|q^2|>|z|>|q^2p|$ and $g_{ij}(z)$ in the region
$|q|>|z|>|qp|$.

This algebra is related
 to the algebra in \cite{FO} \cite{FF}, especially the 
subalgebra of $U_q(\theta,\gtsl_n)$ generated by $x^+_i(z)$, 
although the cubic relations and the Hopf algebra structure were not 
given for their cases. The subalgebras 
generated by $\psi(z)$ and $x^-(z)$ with a comultiplication 
for certain curves were given in \cite{ER}, where they did not
have the full Hopf algebra structure, for example the antipode.
As an algebra,  $U_q(\theta,\gtsl_2)$ is also given in \cite{ER}. 
 
If we take the limit that 
 $p$ goes to zero, we would 
recover $U_q(\hgtsl_n)$. Thus $U_q(\theta,\gtsl_n)$
 gives a deformation of 
$U_q(\hgtsl_n)$. One  question to ask would  concern
functions $h_{ij}(z_1/w, z_2/w)$, which in the case of 
$U_q(\hgtsl_n)$ is $q+q^{-1}$. We note that the functions
$h_{ij}(z_1/w, z_2/w)$ are not constant.
However we maintain this function should have certain geometrical meaning. 
 
Let $\overline{Q}=\oplus_{j=1}^{n-1} \nz \alpha_j$ be the root lattice of 
$\gtsl_n$ and  $\overline{\Lambda}_j=\Lambda_j-\Lambda_0$ be the classical 
part of the $i$-th fundamental weight.

Let Heisenberg algebra be an algebra generated by 
$\{ a_{i,k}|1\leq i \leq n-1,~k\in \nz \setminus \{ 0\} \}$ 
satisfying:
\[ [a_{i,k},a_{j,l}]=
-[b_{i,k},b_{j,l}]=\delta_{k+l,0} 
\frac{[(\alpha_i,\alpha_j)k][k]}{k(1-p^k)}. \]
\[ [a_{i,k}, b_{j,l}]=0\]

Let us  define the  group algebra $\nc(q) [\overline{\cP}]$. 
Let $\overline{\cP}$ be the weight lattice of $\gtsl_n$. 
We fix our free basis $\alpha_2,\cdots, \alpha_{n-1},\overline{\Lambda}_{n-1}$.
They satisfy
\begin{align*}
& e^{\alpha_i}e^{\alpha_j}=(-1)^{(\alpha_i,\alpha_j)}e^{\alpha_j}e^{\alpha_i}
  \qquad 2\leq i,j\leq n-1, \\
& e^{\alpha_i}e^{\overline{\Lambda}_{n-1}}=(-1)^{\delta_{i,n-1}}
  e^{\overline{\Lambda}_{n-1}}e^{\alpha_i}, 
  \qquad 2\leq i\leq n-1.
\end{align*}
For $\alpha=m_2\alpha_2+\cdots m_{n-1}\alpha_{n-1}
     +m_n\overline{\Lambda}_{n-1}$, we set
\[ e^{\alpha}=e^{m_2\alpha_2}\cdots e^{m_{n-1}\alpha_{n-1}}
              e^{m_n\overline{\Lambda}_{n-1}}. \]
Note that the following equations hold.

\begin{align*}
& \overline{\Lambda}_i=-\alpha_{i+1}-2\alpha_{i+2}-\cdots 
                      -(n-i-1)\alpha_{n-1}+(n-i)\overline{\Lambda}_{n-1}, \\
& \alpha_1=-2\alpha_2-3\alpha_3-\cdots -(n-1)\alpha_{n-1}
           +n\overline{\Lambda}_{n-1}.
\end{align*}

Set
\[ \cF_{i}:=\nc(q)[[p^{1/2}]] [a_{j,-k},b_{j,-k}(1\leq j\leq n-1,~k\in \nz_{>0})]\otimes
            \nc(q) [\overline{Q}]
              e^{\overline{\Lambda}_i}
   \quad (0\leq i\leq n-1). \]
This gives the Fock space. 

The action of operators 
$b_{j,k}, a_{j,k},\partial_{\alpha},e^{\alpha}$~
$(1\leq j \leq n-1, \alpha \in \overline{Q})$ is given by 
\begin{align*}
a_{j,k}\cdot f\otimes e^{\beta}& =\begin{cases}
                      a_{j,k}f \otimes e^{\beta}          & k< 0  \\
            \text{$[a_{j,k},f]$}  \otimes e^{\beta} \quad & k> 0
                                  \end{cases}, \\
b_{j,k}\cdot f\otimes e^{\beta}& =\begin{cases}
                      b_{j,k}f \otimes e^{\beta}          & k< 0  \\
            \text{$[b_{j,k},f]$}  \otimes e^{\beta} \quad & k> 0
                                  \end{cases}, \\
\partial_{\alpha}\cdot f\otimes e^{\beta}& 
        =(\alpha,\beta)f\otimes e^{\beta}, 
         \qquad \text{for}~f\otimes e^{\beta}\in \cF_{i} \\
e^{\alpha}\cdot f\otimes e^{\beta} & 
=f\otimes e^{\alpha}e^{\beta}.
\end{align*}

\begin{prop} 
The following gives a highest weight representation on $\cF_{i}$.  
\begin{align*}
\circ & \quad x_j^{\pm}(z)\mapsto
        \exp[\pm \sum_{k>0}\frac{a_{j,-k}q^{\mp \frac{1}{2}k}-
b_{j,-k}(pq^{\pm 1})^{\frac{1}{2}k}} {[k]}z^k]
        \exp[\mp \sum_{k>0}\frac{a_{j,k}q^{\mp \frac{1}{2}k}-
b_{j,k}(pq^{\pm 1})^{\frac{1}{2}k}}
{[k]}z^{-k}] \\
      &  e^{\pm \alpha_j}z^{\pm \partial_{\alpha_j}+1}, \\
\circ & \quad \vphi_j(z)\mapsto
        \exp[-(q-q^{-1})\sum_{k>0}a_{j,-k}z^k]
        \exp[(q-q^{-1})\sum_{k>0}b_{j,k}p^{k/2}z^{-k}]
        q^{-\partial_{\alpha_j}}, \\
\circ & \quad \psi_j(z)\mapsto
        \exp[-(q-q^{-1})\sum_{k>0}b_{j,-k}p^{k/2}z^{k}]
        \exp[(q-q^{-1})\sum_{k>0}a_{j,k}z^{-k}]
        q^{\partial_{\alpha_j}}.
\end{align*}
\end{prop}
\end{ex}

The extension of our construction to other simply laced cases is 
straightforward. 
At this moment, we know very little about these algebras. 
We do not know  when this algebra would  be zero. 
 If we use Drinfeld's
double construction to double the subalgebra generated by 
$x^+_i(z)$ and $\vphi_i(z)$, this will give the whole algebra.
 From the point of view of double construction, an  immediate question is 
 what the universal $R$-matrix is \cite{Dr2}. 
 The answer to this question
will  shed  insight towards understanding this algebra. Another
problem is certainly its representation theory, which hopefully 
should relate to various problem in mathematics and physics. 
In this paper, we  suggest a new  structure inspired 
by the ideas of Drinfeld's current realization of 
affine quantum groups. It is still  a very small step since a complete
understanding of this construction, why it  yields a Hopf algebra 
 and its implication are not yet achieved. We hereby suggest that
a geometric approach towards the realization ofthe representations
of  the algebras should be very helpful.  \cite{GRV}\cite{Go} 

\begin{ack}
J.D. would like to thank I.B. Frenkel for his advice and encouragement. 
We would also like to thank T. Miwa and B. Feigin for useful disscussions. 
J.D. is supported by the grant Reward research (A) 08740020 from the 
Ministry of Education of Japan. K.I. is partly supported by the Japan Society
for the Promotion of Science. 
\end{ack}

\end{document}